\documentclass[11pt]{article}
\usepackage{amssymb}

\newcommand{\Z}{{\mathbb Z}}

\newcommand{\R}{{\mathbb R}}

\newcommand{\SUM}{\raisebox{-0.4ex}{\mbox{\Large $\Sigma$}}}
\newcommand{\eps}{\varepsilon}

\newtheorem{theorem}{Theorem}

\newtheorem{proposition}{Proposition}

\title{An application of linear programming duality to discrete 
Fourier analysis and additive problems}

\author{Ernie Croot}

\begin{document}

\maketitle

\section{Introduction}

Suppose that $p$ is a prime number and that 
$$
S\ \subseteq\ \Z_p.
$$
Associate to $S$ the indicator function $S(n)$, which equals $1$ if 
$n \in S$ and equals $0$ if $n \not \in S$.

As is well known, the additive properties of $S$ are strongly dependent
on the Fourier coefficients
$$
\hat S(a)\ :=\ \SUM_{n \in \Z_p} S(n) \omega^{an}\ =\ \SUM_{n \in S} \omega^{an},
\ {\rm where\ } \omega = e^{2\pi i/p};
$$
in particular, if the size of the second-largest
Fourier coefficient is ``small'', then $|S+S|$ must be appreciably larger than $|S|$.  

Now let us suppose that we want to show that a sumset $S+S$ is large, or perhaps
we wish to prove that $S$ has some other additive property, such as that it
contains many three-term arithmetic progressions.  Although $S$ itself may have
a ``large'' second-largest Fourier coefficient, we can
imagine that perhaps there exists a ``large'' subset
$$
T\ \subseteq\ S
$$
such that the second-largest Fourier coefficient of the convolution $S*T$ is ``small''.
If so, then one can show that (for appropriate notions of ``large'' and ``small'') this
implies that $S+T$ is ``large'', and therefore so is
$$
|S+S|\ \geq\ |S+T|.
$$
Recall that the convolution $S*T$ defined by
$$
(S*T)(n)\ :=\ \SUM_{a+b=n} S(a)T(b)
$$
has the properties 
$$
(S*T)(n)\ >\ 0\ \ \iff\ \ n \in S+T,
$$
and
$$
\widehat{(S*T)}(a)\ =\ \hat S(a) \hat T(a).
$$
\bigskip

Actually, in place of the set function $T(n)$, all we really need to do is to produce a
function
$$
f\ :\ \Z_p\ \to\ \R,
$$
where
\begin{eqnarray} \label{f1}
n \in S\ &\Longrightarrow&\ f(n) \geq 0;\ {\rm and,} \nonumber \\ 
n \in S^c\ &\Longrightarrow&\ f(n) \leq 0,
\end{eqnarray}
because for such $f$ we will have that
$$
(f*S)(n)\ >\ 0\ \Longrightarrow\ n \in S+S.
$$
If we in addition had that 
\begin{equation} \label{f2}
\hat f(a)\ =\ 0\ \ {\rm at\ all\ places\ }a\ {\rm where\ } |\hat S(a)|\ {\rm is\ ``large"},
\end{equation}
then we would have some control over the size of the largest non-zero
Fourier coefficient of $f*S$.  
\bigskip

It reasonable to expect that in a lot of instances we can produce 
a function $f$ satisfying (\ref{f1}) and (\ref{f2}), provided that 
there aren't too many places $a$ where
$|\hat S(a)|$ is ``large''.  However, what is not so obvious is that, even 
when we demand that $\hat f(a) = 0$ at a lot of places $a$, if no such
function $f$ exists, we still can get a rather nice and useful conclusion 
by applying the principle of the seperating hyperplane, which is a basic type of 
duality principle from Linear Programming used to prove Farkas's Lemma.   
In order to state informally what  our result gives, let us introduce the 
following definition:
\bigskip

\noindent {\bf Definition.}  We say that a function $f : \Z_p \to \R$ is a 
{\it generalized balanced function} for some set $S$ if it satisfies the following
properties:

$\bullet$ if $x \in S$, then $f(x) \geq 0$;

$\bullet$ if $x \in S^c := \Z_p \setminus S$, then $f(x) \leq 0$; 

$\bullet$ $\SUM_x f(x)\ =\ 0$.
\bigskip

\noindent Our theorem  will say that given a sequence of places 
$a_1,...,a_k \neq 0$, and given some function $g : \Z_p \to [0,1]$, we can
either find a non-zero {\it generalized balanced function} for ${\rm support}(g)$
whose Fourier transform vanishes at $a_1,...,a_k$ (and satisfies some
additional constraints); or, we can find a 
function that behaves like a generalized balanced function for 
${\rm support}(g)$ in the sense that it satisfies the first two bullets above for
``most'' $x \in \Z_p$, and has the additional, very useful property that
its Fourier transform has small support.  Applications of this theorem to 
additive number theoretical problems -- in particular, analyzing the additive
properties of level sets of sumsets, where by ``level set'' we mean something
like
$$
\{ n \in \Z_p\ :\ (S*S)(n)\ \in\ [L-\eps,\ L + \eps]\}
$$
or perhaps a triple convolution $S*S*S$ -- will perhaps be worked out in a 
forthcoming paper.

\begin{theorem} \label{main_theorem} 
Suppose that 
$$
g\ :\ \Z_p\ \to\ [0,1],
$$
and suppose that 
$$
a_1, a_2, ..., a_k\ \in\ \Z_p \setminus \{0\}
$$
are any $k$ distinguished non-zero places.  Let $E$ be some integer satisfying
$$
E\ \geq\ 0.
$$
Then, one or the other of the following two conclusions must hold:
\bigskip

$\bullet$  {\bf (Vanishing Generalized Balanced Function)} Either there exists a function 
$$
h\ :\ \Z_p\ \to\ [-1,1],
$$
satisfying
\begin{eqnarray*}
n \in {\rm support}(g)\ &\Longrightarrow&\ h(n)\ \geq\ 0, \\
n \in {\rm support}(g)^c\ &\Longrightarrow&\ h(n)\ \leq\ 0,
\end{eqnarray*}
and
$$
\SUM_n h(n)\ =\ 0;\ ||h||_1\ =\ \SUM_n |h(n)|\ \geq\ E;
$$
and
$$
\hat h(a_1)\ =\ \hat h(a_2)\ =\ \cdots\ =\ \hat h(a_k)\ =\ 0.
$$
\bigskip

$\bullet$ {\bf (Generalized Balanced Function with Small Spectral Support)} 
Or, there exists a function 
$$
h\ :\ \Z_p\ \to\ \R,
$$
such that 
$$
{\rm support}(\hat h)\ \subseteq\ \{0\} \cup \{a_1,...,a_k\} \cup \{-a_1,...,-a_k\},
$$
and such that, apart from at most 
$$
(2k+1)E\ \ {\rm exceptions},
$$
we will have that
\begin{eqnarray*}
n \in {\rm support}(g)\ &\Longrightarrow&\ h(n)\ >\ 0,\ {\rm and} \\
n \in {\rm support}(g)^c\ &\Longrightarrow&\ h(n)\ <\ 0.
\end{eqnarray*}
\end{theorem}

\noindent {\bf Remark 1.}  We note that when $E = 0$ the Theorem is trivially
true, since the $0$ function satisfies the first conclusion in that case;
also, when $E \geq (2k+1)^{-1} p$ the second conclusion is trivially true. 
\bigskip

\noindent {\bf Remark 2.}  It would be nice to have a theorem where in place
of the first conclusion above we had one where we have some control over
the sum over $n$ of $h(n)$, such as
$$
\SUM_n h(n)\ >\ F,
$$
for some function $F$ that depends on, say, 
$$
|\hat g(a_1)|,\ ...,\ |\hat g(a_k)|.
$$
It might be possible to prove a theorem like this by developing some quantitative
version of Farkas's Lemma, and using the method of proof in the present paper.
This would undoubtedly have many nice applications, and could possibly lead
to a new proof of Roth's Theorem on three-term arithmetic progressions.

\section{Proof of Theorem \ref{main_theorem}} \label{proof_section}

\subsection{Separating hyperplanes}

As mentioned earlier, we will require the following basic proposition.  Its proof requires the 
{\it principle of the seperating hyperplane}, which implies that if a convex hull $H$
of a some points in $\R^m$ does not contain some point $P$, then there exists a 
hyperplane that separates $\R^m$ into three regions:  One region contains $P$,
another region contains $H$, and the third region is the hyperplane itself.

\begin{proposition} \label{simplex_corollary}
Suppose that $M$ is an $m \times n$ matrix with real entries,
where $n > m$.  Then, one of the following must hold:
\bigskip

$\bullet$ Either there exists a non-negative vector 
$$
v\ =\ (v(1),...,v(n))\ \in\ \R_{\geq 0}^n\ \ {\rm with\ \ }\ \SUM_j v(j)\ =\ 1,
$$ 
having at most $m$ non-zero entries, and satisfying
$$
Mv\ =\ 0;
$$ 
\bigskip

$\bullet$ or, there exists a vector 
$$
w\ \in\ \R^m,
$$
such that 
$$
w M\ \in\ \R_{> 0}^n.
$$ 
\end{proposition}

\noindent {\bf Proof of the Proposition.}  Consider the convex hull $H$ of the 
columns of $M$.  Every point of this convex hull is a linear combination of these
columns, where the coefficients are all $\geq 0$ and sum to $1$.   
There are two possibilities:  Either $0 \in H$ or $0 \notin H$.  

First suppose that $0 \in H$.  Then, by taking a simplicial decomposition 
of $H$, we find that there exists a simplex consisting of at most $m$ vertices 
drawn from the column vectors of $M$, which contains $0$.  To say that $0$
lies in or on this simplex means that some linear combination of the
these $\leq m$ vertex vectors, using non-negative coefficients that sum to $1$,
sums to $0$.  Expressing this in matrix and vector notation, we obtain the
first conclusion of the Proposition.

Now suppose that $0 \notin H$.  By the principle of the separating hyperplane, there
exists a hyperplane of $\R^m$ such that $0$ is on one side of the hyperplane,
while $H$ is on the other.  Let $w$ be a normal vector to this hyperplane so 
that if $x \in H$ then 
$$
w \cdot x\ >\ 0.
$$
It follows that
$$
w M\ \in\ \R_{> 0}^n.
$$

\subsection{Body of the proof of Theorem \ref{main_theorem}}

We now prove Theorem \ref{main_theorem} by applying the above 
proposition iteratively.  First, let 
$$
\{b_1,...,b_t\}\ :=\ \{a_1,...,a_k\}\ \cap\ \{1,...,(p-1)/2\}.
$$
If we can produce a function $h : \Z_p \to \R$ such that $\hat h$ vanishes at these 
places $b_i$, then it will automatically at the negatives of these places, from the fact that
$$
\hat h(a)\ =\ \overline{\hat h(-a)}.
$$

We will construct a sequence of matrices 
$$
M_1,\ M_2,\ ...,\ M_T,
$$
where 
$$
M_i\ =\ m \times n_i,
$$
and a sequence of vectors
$$
v_1,\ v_2,\ ...,\ v_T\ \in\ \R_{\geq 0}^n
$$
(or maybe the last vector is only $v_{T-1}$), and 
then we will read off properties of $M_T$ and $v_T$ to prove our theorem.
\bigskip

We begin by letting 
$$
m\ :=\ 2t+1,\ n\ :=\ p,
$$
and then we define $M_1$ to be the $m \times n$ matrix whose $j$th column 
is given as follows:  First, if $j \in {\rm support}(g)$, then the column vector is
\begin{eqnarray*}
(1,\ \cos(2\pi j b_1/p),&& \sin(2\pi j b_1/p),\ \cos(2\pi j b_2/p),\ \sin(2\pi j b_2/p),\\
&& ...,\ \cos(2\pi j b_t/p),\ \sin(2\pi j b_t/p)),
\end{eqnarray*}
and if $j \not \in {\rm support}(g)$, then the column vector is
\begin{eqnarray*}
(-1,\ -\cos(2\pi j b_1/p), && -\sin(2\pi j b_1/p),\ -\cos(2\pi j b_2/p),\ -\sin(2\pi j b_2/p),\ ...,\\ 
&& -\cos(2\pi j b_t/p), - \sin(2\pi j b_t/p) ).
\end{eqnarray*}
\bigskip

Given that we have constructed $M_r$, and that our iterative process (described below) 
did not end with $M_r$, we apply Proposition \ref{simplex_corollary}
with $M := M_r$.  So, one or the other of the conclusions of that Proposition must hold.

\subsubsection{Case 1 (first conclusion of Proposition holds)}
  
Let us first suppose that the first conclusion of the Proposition holds, and let $v$ be 
the vector appearing there.  From $v$, which has $n_r$ coordinates,
we produce a vector $v_r$ having $n$ coordinates as follows:  First, the columns
of $M_r$ correspond to particular columns of $M$, and let us say that the $j$th
column of $M_r$ corresponds to the $c_j$th column of $M$.  Then, writing 
$$
v\ =\ (v(1),\ ...,\ v(n_r)),
$$
we define
$$
v_r(c_j)\ =\ v(j);\ {\rm and,\ for\ } i \not \in \{c_1,...,c_{n_r}\},\ {\rm we\ set\ } 
v_r(i)\ =\ 0.
$$
So, basically the coordinates of $v_r$ that correspond to columns that were deleted
when passing from $M_1$ to $M_r$ are set to $0$, while the coordinates corresponding
to the other columns of $M_1$  (that were not deleted) are assigned their respective
values from the vector $v$. 

We note that 
\begin{equation} \label{Mvr}
M_1 v_r\ =\ M_r v\ =\ 0.
\end{equation}
\bigskip

If 
$$
r\ <\ E,
$$
then we define the matrix $M_{r+1}$ by taking $M_r$ and removing the 
columns corresponding to places where $v$ has a non-zero entry.  Furthermore,
we let $n_{r+1}$ be the number of columns of $M_{r+1}$. 

On the other hand, if 
$$
r\ =\ E,
$$
then we STOP the process of generating matrices $M_j$ and vectors $v_j$, and
set $T := r$.  We note that the non-zero coordinates of the vectors 
$$
v_1,...,v_T
$$ 
are all mutually disjoint, and so letting 
$$
V\ :=\ v_1 + \cdots + v_T
$$
we will have from (\ref{Mvr}) that 
$$
V \in \R_{\geq 0}^n,\ M_1 V\ =\ 0,\ ||V||_\infty\ =\ 1,\ {\rm and\ } 
||V||_1\ \geq\ T\ =\ E.
$$

So, if 
$$
V\ =\ (V(1),\ V(2),\ ...,\ V(p)),
$$
then if we define
$$
h(a)\ =\ \left \{\begin{array}{rl} V(a),\ &{\rm if\ } a\in {\rm support}(g); \\
-V(a),\ & {\rm if\ } a\notin {\rm support}(g), \end{array}\right.
$$
we will have that since, again, the supports of the $v_i$ are all disjoint,
$$
||h||_1\ =\ \SUM_{i=1}^T ||v_i||_1\ \geq\ E,
$$
and
\begin{eqnarray*}
a\ \in\ {\rm support}(g)\ &\Longrightarrow&\ h(a)\ \geq\ 0;\ {\rm and,} \\
a\ \in\ {\rm support}(g)^c\ &\Longrightarrow&\ h(a)\ \leq\ 0.
\end{eqnarray*}

Furthermore, although it takes a little work to see, one can read off from the
fact that 
$$
M_1 v_1\ =\ M_1 v_2\ =\ \cdots\ =\ M_1 v_T\ =\ 0,
$$
the conclusions
$$
\hat h(a_1)\ =\ \cdots\ =\ \hat h(a_k)\ =\ 0,\ {\rm and\ } \SUM_a h(a)\ =\ 0.
$$
This then would give the first conclusion claimed by our Theorem.

\subsubsection{Case 2 (second conclusion of Proposition holds)}

If the second conclusion of the Proposition holds, then there exists a vector
$$
w\ \in\ \R^m,\ w\ \neq\ 0,
$$
such that 
\begin{equation} \label{wM}
w M_r\ \in\ \R_{> 0}^{n_r}.
\end{equation}
Letting 
$$
w\ =\ (w(1),\ w(2),\ ...,\ w(m)),
$$
we find that (\ref{wM}) is equivalent to the following:  First, if 
$$
x\ \in\ {\rm support}(g),
$$ 
and $x$ does not correspond to one of the columns that was deleted in passing
from $M_1$ to $M_r$, then 
\begin{equation} \label{inv_fourier1}
w(1)\ +\ \SUM_{j=1}^t \left ( w(2j) \cos(2\pi x b_j/p) + w(2j + 1) \sin(2\pi x b_j/p) \right )\ 
\geq\ 0.
\end{equation}
And second, if 
$$
x\ \in\ {\rm support}(g)^c,
$$
and, again, $x$ does not correspond to a deleted column, then
\begin{equation} \label{inv_fourier2}
w(1)\ +\ \SUM_{j=1}^t (w(2j) \left (\cos(2\pi x b_j/p) + w(2j + 1)\sin(2\pi x b_j/p) \right )\ 
\leq\ 0.
\end{equation}

We can think of (\ref{inv_fourier1}) and (\ref{inv_fourier2}) as the inverse Fourier
transform of a certain function 
$$
h\ :\ \Z_p\ \to\ \R,
$$
where
$$
\hat h(a)\ =\ \left \{ \begin{array}{rl} p w(1),\ & {\rm if\ } a = 0; \\
p w(2j)/2 - i p w(2j+1)/2,\ & {\rm if\ } a=b_j; \\
p w(2j)/2 + i p w(2j+1)/2,\ & {\rm if\ } a = -b_j,
\end{array}\right.
$$

One can check that
\begin{eqnarray*}
h(x)\ &=&\ w(1)\ +\ \SUM_{j=1}^t (w(2j)/2 - i w(2j+1)/2) e^{2\pi i x b_j/p}\\
&&\hskip0.5in +\ \SUM_{j=1}^t (w(2j)/2 + i w(2j+1)/2) e^{-2\pi i x b_j/p},
\end{eqnarray*}
which, along with (\ref{inv_fourier1}) and (\ref{inv_fourier2}), implies that,
apart from at most $Tm$ exceptions (which is the maximum number of deleted
columns), we will have that 
\begin{eqnarray*}
x \in {\rm support}(g)\ &\Longrightarrow&\ h(x)\ >\ 0;\ {\rm and,} \\
x \in {\rm support}(g)^c\ &\Longrightarrow&\ h(x)\ <\ 0.
\end{eqnarray*}
This then finishes the proof of our theorem.

\end{document}